\documentclass[12pt]{article} \usepackage{amsmath, amssymb,latexsym, epsfig, color}

\newtheorem{theorem}{Theorem}[section]
\newtheorem{proposition}[theorem]{Proposition}

\newtheorem{definition}[theorem]{Definition}
\newtheorem{conjecture}[theorem]{Conjecture}
\newtheorem{remark}[theorem]{Remark} \newtheorem{lemma}[theorem]{Lemma}
\newtheorem{example}[theorem]{Example} 

\newenvironment{numtheorem}[1]{\medskip\noindent {\bf Theorem~#1.}\em}{}

\def\proof{\smallskip\noindent {\it Proof: \ }} \def\endproof{\hfill$\square$\medskip}

  \newcommand{\lk}{\mbox{\upshape lk}\,}     \newcommand{\field}{{\bf k}}  \newcommand{\Skel}{\mbox{\upshape Skel}\,}          \newcommand{\MON}{\mbox{\upshape \MON}}     \newcommand{\codim}{\mbox{\upshape codim}\,} \newcommand{\lcm}{\mbox{\upshape lcm}\,} \newcommand{\mubc}{multiplicity upper bound conjecture} \newcommand{\mlbc}{multiplicity lower bound conjecture} \title{Face ring multiplicity via CM-connectivity sequences}

\author{Isabella Novik \thanks{Research partially supported by NSF grants
DMS-0500748 and SBE-0123552}\\ \small Department of Mathematics, Box 354350\\[-0.8ex] \small University of Washington, Seattle, WA 98195-4350, USA,\\[-0.8ex] \small \texttt{novik@math.washington.edu} \and Ed Swartz \thanks{Research partially supported by NSF grant DMS-0245623}\\ \small Department of Mathematics, \\[-0.8ex] \small Cornell University, Ithaca NY, 14853-4201, USA, \\[-0.8ex] \small \texttt{ebs@math.cornell.edu } }

\begin{document}

\maketitle \begin{abstract} The multiplicity conjecture of Herzog, Huneke, and Srinivasan is verified for the face rings of the following classes of simplicial complexes: matroid complexes, complexes of dimension one and two, and Gorenstein complexes of dimension at most four. The lower bound part of this conjecture is also established for the face rings of all doubly Cohen-Macaulay complexes whose 1-skeleton's connectivity does not exceed the codimension plus one as well as for all $(d-1)$-dimensional $d$-Cohen-Macaulay complexes. The main ingredient of the proofs is a new interpretation of the minimal shifts in the resolution of the face ring $\field[\Delta]$ via the Cohen-Macaulay connectivity of the skeletons of $\Delta$.
\end{abstract}

\section{Introduction}
In this paper we prove the multiplicity conjecture of Herzog, Huneke, and Srinivasan for the face rings of several classes of simplicial complexes.

Throughout the paper we work with the polynomial ring
   $S=\field[x_1, \ldots, x_n]$
over an arbitrary field $\field$. If $I\subset S$ is a homogeneous ideal, then the (bi-graded) {\bf Betti numbers} of $I$, $\beta_{i,j}=\beta_{i,j}(I)$,
   are the invariants that appear
in the minimal free resolution of $S/I$ as an $S$-module:
\[ 0 \rightarrow
\bigoplus_j S(-j)^{\beta_{l,j}(I)} \rightarrow \ldots \rightarrow \bigoplus_j S(-j)^{\beta_{2,j}(I)} \rightarrow \bigoplus_j S(-j)^{\beta_{1,j}(I)} \rightarrow S \rightarrow S/I \rightarrow 0 \] Here $S(-j)$ denotes $S$ with grading shifted by $j$ and $l$ denotes the length of the resolution. In particular, $l \geq \codim(I)$.

Our main objects of study are the {\bf maximal and minimal shifts} in the resolution of $S/I$ defined by $M_i=M_i(I)=\max\{j : \beta_{i,j}\neq 0\}$ and $m_i=m_i(I)=\min \{j : \beta_{i,j}\neq 0\}$ for $i=1, \ldots, l$, respectively.
The following conjecture due to Herzog, Huneke, and Srinivasan \cite{HerzSr98} is known as the multiplicity conjecture.

\begin{conjecture}  \label{multiplicity-conj} Let $I\subset S$ be a homogeneous ideal of codimension $c$. Then the multiplicity of $S/I$, $e(S/I)$, satisfies the following upper bound:
$$ e(S/I) \leq ( \prod_{i=1}^{c} M_i )/c!.
$$
Moreover, if $S/I$ is Cohen-Macaulay, then also

$$ e(S/I) \geq ( \prod_{i=1}^{c} m_i )/c!.
$$
\end{conjecture}

This conjecture was motivated by the result due to Huneke and Miller \cite{HuMill} that if $S/I$ is Cohen-Macaulay and $M_i=m_i$ for all $i$ (in such a case $S/I$ is said to have a {\bf pure resolution}), then $e(S/I)= (\prod_{i=1}^{c} m_i )/c!$. Starting with the  paper of Herzog and Srinivasan \cite{HerzSr98} a tremendous amount of effort was put in establishing Conjecture \ref{multiplicity-conj} for various classes of rings $S/I$ (see \cite{Franc, Gold03, GoldSchSr,  GuardoTuyl,  HerzSr98, HerzZheng, MiglNagRom,  Miro,  Rom05, Sr98} and the survey article \cite{FrancSr}).
In particular, the conjecture was proved in the following cases: $S/I$ has a {\bf quasi-pure resolution}  (that is, $m_i(I) \geq M_{i-1}(I)$ for all $i$) \cite{HerzSr98}, $I$ is a stable or squarefree strongly stable ideal \cite{HerzSr98}, $I$ is a codimension 2 ideal \cite{Gold03, HerzSr98, Rom05}, and $I$ is a codimension 3 Gorenstein ideal \cite{MiglNagRom}.

We investigate Conjecture \ref{multiplicity-conj} for squarefree monomial ideals or, equivalently, face ideals of simplicial complexes.
If $\Delta$ is a simplicial complex on the vertex set $[n]=\{1, 2, \ldots, n\}$, then its {\bf face ideal} (or the {\bf Stanley-Reisner ideal}), $I_\Delta$, is the ideal generated by the squarefree monomials corresponding to non-faces of $\Delta$, that is, $$ I_\Delta = \langle
       x_{i_1}\cdots x_{i_k} \, : \, \{i_1<\cdots <i_k\}\notin \Delta
              \rangle, $$
and the {\bf face ring} (or the {\bf Stanley-Reisner ring}) of $\Delta$ is $\field[\Delta]:=S/I_\Delta$ \cite{St96}.

Various combinatorial and topological invariants of $\Delta$ are encoded in the algebraic invariants of $I_\Delta$ and vice versa \cite{BrHerz, St96}.
The Krull dimension of $\field[\Delta]$, $\dim  \field[\Delta]$, and the topological dimension of $\Delta$, $\dim \Delta$, are related by
$\dim  \field[\Delta] =\dim \Delta +1$ and so $$       \codim(I_\Delta)=n-\dim
\Delta-1.
$$
The Hilbert series of $\field[\Delta]$ is determined by knowing the number of faces in each dimension.  Specifically, let  $f_i$  be the number of $i$-dimensional faces. By convention, the empty set is the unique face of dimension minus one.  Then, $$\sum^\infty_{i=0} \dim_\field \field[\Delta]_i \lambda^i = \frac{h_0 + h_1 \lambda +\dots + h_d \lambda^d}{(1-\lambda)^d},$$ where, $\field[\Delta]_i$ is the $i$-th graded component of $\field[\Delta]$, $d=\dim\Delta+1=\dim  \field[\Delta]$, and
   \begin{equation} \label{h-vector}
    h_i = \sum^i_{j=0} (-1)^{i-j} \binom{d-j}{d-i} f_{j-1}.
\end{equation}
The multiplicity $e(\field[\Delta])$  equals the number of top-dimensional faces of $\Delta$ which in turn is $h_0 + \dots + h_d.$ The minimal and maximal shifts have the following interpretation in terms of the reduced homology:
\begin{eqnarray}  \label{M-interpr}
M_i(I_\Delta)&=&\max\{|W| \, : \, W\subseteq[n] \mbox{ and }
              \tilde{H}_{|W|-i-1}(\Delta_W; \field)\neq 0\}, \\ \label{m-interpr} m_i(I_\Delta)&=&\min\{|W| \, : \, W\subseteq[n] \mbox{ and }
              \tilde{H}_{|W|-i-1}(\Delta_W; \field)\neq 0\}.
\end{eqnarray}
(Here $\Delta_W$ denotes the {\bf induced subcomplex} of $\Delta$ whose vertex set is $W$, and $\tilde{H}(\Delta_W; \field)$ stands for the reduced simplicial homology of $\Delta_W$ with coefficients in $\field$.
The above expressions for $M_i$ and $m_i$ follow easily from Hochster's formula on the Betti numbers $\beta_{i,j}(I_\Delta)$ \cite[Theorem II.4.8]{St96}.) Thus, for  face ideals, Conjecture \ref{multiplicity-conj} can be considered as a purely combinatorial-topological statement. As we will see below, the upper bound part of the conjecture is closely related to the celebrated Upper Bound Theorem for polytopes and Gorenstein* complexes \cite{St75}.

In this paper we prove Conjecture \ref{multiplicity-conj} for the face rings of the following classes  of simplicial complexes.
\begin{itemize}
\item Matroid complexes. A simplicial complex is called a {\bf matroid complex} if it is {\bf pure} (that is, all its maximal under inclusion faces have the same dimension) and all  its induced subcomplexes are pure.

\item One- and two-dimensional complexes.

\item Three- and four-dimensional Gorenstein complexes. \end{itemize} The first result about matroid complexes has a flavor similar to that of squarefree strongly stable ideals, while the last two results complement the fact that the multiplicity conjecture holds for codimension 2 ideals and for codimension 3 Gorenstein ideals.

Recall that a simplicial complex is
called {\bf Gorenstein} over $\field$ if its face ring is Gorenstein.
Similarly, a simplicial complex $\Delta$ is said to be {\bf Cohen-Macaulay over $\field$} (CM, for short) if its face ring $\field[\Delta]$ is Cohen-Macaulay.
A simplicial complex is {\bf $q$-Cohen-Macaulay} ($q$-CM, for short) if for every set $U\subset [n]$, $0\leq |U|\leq q-1$, the induced subcomplex $\Delta_{[n]-U}$ (that is, the complex obtained from $\Delta$ by removing all vertices in $U$) is a CM complex of the same dimension as $\Delta$. 2-CM complexes are also called {\bf doubly CM complexes}.

By Reisner's criterion \cite{Reisner}, $\Delta$ is CM  if and only if $$
   \tilde{H}_i(\lk F; \field) = 0 \qquad \mbox{for all } F\in\Delta \mbox{ and }
              i<\dim\Delta -|F|.
$$
Here $\lk F=\{G\in \Delta \, : \, F\cap G=\emptyset, \, F\cup G \in \Delta\}$ is the {\bf link of face $F$ in $\Delta$} (e.g., $\lk \emptyset=\Delta$). Thus, a 1-dimensional complex (that is, a graph) is CM if and only if it is connected.
Moreover, a 1-dimensional complex is $q$-CM if it is $q$-connected in the usual graph-theoretic sense. A 0-dimensional complex on $n$ vertices is $n$-CM.

We verify the lower-bound part of Conjecture \ref{multiplicity-conj} for the face rings of the following classes of simplicial complexes \begin{itemize} \item 2-CM $(d-1)$-dimensional simplicial complexes on $n$ vertices whose 1-skeleton is at most $(n-d+1)$-connected.
\item d-CM $(d-1)$-dimensional complexes.
\end{itemize}
The restriction that the 1-skeleton is at most $(n-d+1)$-connected is a rather mild one: it means that there exists a subset $W\subset [n]$ of size $d-1$ such that $\Delta_W$ is disconnected.
This condition is satisfied, for instance, by the order complex of an arbitrary graded poset of rank $d$ with $\geq d^2$ elements.

The key notion in the proof of the lower bound part of Conjecture \ref{multiplicity-conj} is that of the
   {\bf CM-connectivity sequence of a CM complex}.
   Recall that the $i$-skeleton of a simplicial complex $\Delta$, $\Skel_i(\Delta)$,
   is the collection of all faces of $\Delta$ of dimension $\leq i$.

\begin{definition}    \label{CM-connectivity-seq}
For a CM simplicial complex $\Delta$ define $q(\Delta):=\max \{q \, : \, \Delta \mbox{ is $q$-CM} \}$, and $q_i=q_i(\Delta):=q(\Skel_i(\Delta))$ for $0\leq i \leq \dim\Delta$.
The sequence
$
(q_0, \ldots, q_{\dim(\Delta)})
$
is called the {\bf CM-connectivity sequence} of $\Delta$. \end{definition}

A recent result of Fl{\o}ystad \cite[Corollary 2.4]{Floyst} implies that the CM-connectivity sequence is strictly decreasing. Thus for a $(d-1)$-dimensional CM complex $\Delta$ on $n$ vertices, $$1 \leq q_{d-1} < q_{d-2} < \cdots < q_1 <q_0 =n.$$ On the other hand, it follows easily from the minimality of the resolution (see \cite[Proposition 1.1]{BrunsHibi}) that for any  $(d-1)$-dimensional complex $\Delta$ on $n$ vertices $$
2 \leq m_1(I_\Delta) < m_2(I_\Delta) \cdots < m_{n-d}(I_\Delta) \leq n.
$$
This suggests that there may be a connection between the CM-connectivity sequence and the sequence of the minimal shifts.
Existence of such a connection is one of the two main ingredients of our proofs.

\begin{theorem}   \label{main_thm}
Let $\Delta$ be a  CM complex on $n$ vertices, and let $(q_0, \ldots, q_{d-1})$ be the CM-connectivity sequence of $\Delta$, where $d-1$ is the dimension of $\Delta$.
Then $$[n]-\{m_1(I_\Delta), \ldots, m_{n-d}(I_\Delta) \} = \{n-q_0+1, n-q_1+1, \ldots, n-q_{d-1}+1\}, $$ and hence $$ \frac{\prod_{i=1}^{n-d} m_i(I_\Delta)}{(n-d)!} =\frac{n!}{(n-d)!} \cdot \frac{\prod_{i=1}^{n-d} m_i(I_\Delta)}{n!}= \frac{n(n-1)\cdots(n-(d-1))}{(n-q_0+1)\cdots (n-q_{d-1}+1)}. $$ \end{theorem} (In view of this theorem, we refer to $n-q_i+1$ as the {\bf $i$-th skip} in the
$m$-sequence.)

The second ingredient is the following purely combinatorial fact.

\begin{theorem}  \label{q-estimate}
Let $\Delta$ be a  CM $(d-1)$-dimensional complex, and let $(q_0, \ldots, q_{d-1})$ be its CM-connectivity sequence.
Then
$$
e(\field[\Delta])=f_{d-1}(\Delta) \geq \frac{q_0q_1 \cdots q_{d-1}}{d!}.
$$
\end{theorem}

Combining the last two theorems we obtain that the lower bound part of Conjecture \ref{multiplicity-conj} holds for all CM complexes whose CM-connectivity sequence satisfies $\frac{n(n-1)\cdots(n-(d-1))}{(n-q_0+1)\cdots (n-q_{d-1}+1)} \leq \frac{q_0q_1 \cdots q_{d-1}}{d!}$. We then work out which sequences satisfy this inequality.

When $\field[\Delta]$ has a pure resolution, the combination of Huneke and
Miller's formula for rings with pure resolutions and Theorems  \ref{main_thm}, and
\ref{q-estimate} puts a strong restriction on the  CM-connectivity sequence
of $\Delta.$  Indeed,  $(q_0,\dots,q_{d-1})$ must satisfy

$$\frac{q_0 \cdots q_{d-1}}{d!} \le \frac{n (n-1) \cdots (n-d+1)}{(n-q_0+1)
\cdots (n-q_{d-1}+1)}.$$
For instance, when $d=2,$ it is immediate that the only possible values for
$q_1$ are $1,2$ or $n-1.$  Each of these values does occur as trees,
circuits and  complete graphs all have pure resolutions \cite{BrunsHibi98}.

It is worth mentioning that there is a more recent conjecture \cite{HerzZheng, MiglNagRom2} asserting that if for a ring $S/I$, the multiplicity of $S/I$ equals the lower bound or the upper bound of Conjecture \ref{multiplicity-conj}, then $S/I$ is Cohen-Macaulay and has a pure resolution. For the case of Cohen-Macaulay rings with a quasi-pure resolution this conjecture was established in \cite[Theorem 3.4]{HerzZheng}.
Here we verify this conjecture in all the cases where we are able to prove the original multiplicity conjecture (Conjecture \ref{multiplicity-conj}) with the exception of the upper bound for 2-dimensional complexes.

The structure of the paper is as follows. Section 2 is devoted to the class of matroid complexes.
   In Section 3 we verify
Theorems \ref{main_thm} and \ref{q-estimate} as well as derive their application to 2-CM complexes with $q_1\leq n-d+1$ and $d$-CM complexes. In Section 4 we treat one- and two-dimensional complexes. Finally in Section 5 we discuss three- and four-dimensional Gorenstein complexes.

\section{Matroid complexes}            \label{matroid_section}
In this section we verify Conjecture \ref{multiplicity-conj} (together with the treatment of equality) for matroid complexes, namely we prove the following theorem. In the rest of the paper we abuse notation and write $m_i(\Delta)$ and $M_i(\Delta)$ instead of $m(I_\Delta)$ and $M_i(I_\Delta)$, respectively.

\begin{theorem}  \label{matroid-multiplicity-thm} Let $\Delta$ be a $(d-1)$-dimensional matroid complex on $n$ vertices.
Then
$$
(\prod_{i=1}^{n-d} m_i(\Delta))/(n-d)! \leq  f_{d-1}(\Delta) \leq (\prod_{i=1}^{n-d} M_i(\Delta))/(n-d)!.
$$
Moreover, if one of the bounds is achieved, then $\field[\Delta]$ has a pure resolution.
\end{theorem}

We start by reviewing necessary background on matroid complexes.
A complex $\Delta$ on the vertex set $[n]$ is a matroid complex if for every subset $\emptyset \subseteq W \subseteq [n]$, the induced subcomplex $\Delta_W$ is pure.
Equivalently (see \cite[Proposition~2.2.1]{NicolettiWhite}), a matroid complex is a complex that consists of the independent sets of a matroid.

The following lemma summarizes several well known properties of matroid complexes.

\begin{lemma}   \label{matroid-properties}
Let $\Delta$ be a matroid complex. Then
\begin{enumerate}
\item Every induced subcomplex of $\Delta$ is a matroid complex.
\item $\Delta$ is CM. Moreover,  $\Delta$ is  the join of a simplex and a
        2-CM complex.
\item $\Delta$ is either a cone or has a non-vanishing top homology.
\end{enumerate}
\end{lemma}
\proof Part 1 is obvious from the definition of a matroid complex.
For Part 2 see  Proposition III.3.1 and page 94 of \cite{St96}.
Part 3 is a consequence of Part 2 and the fact that a 2-CM complex has a nonvanishing top homology (see top of page 95 in \cite{St96}).
\endproof

Another fact we need for the proof of Theorem \ref{matroid-multiplicity-thm} is \begin{lemma}  \label{M-and-m-properties} Let $\Delta$ be a $(d-1)$-dimensional complex on the vertex set $[n]$. Then \begin{enumerate} \item $M_{n-d}(\Delta)=n$ unless $\tilde{H}_{d-1}(\Delta; \field)=0$, in which case $M_{n-d}(\Delta)<n$.
\item $M_i(\Delta)=\max \{ M_i(\Delta_{[n]-x}) : x\in[n] \}$   and
$m_i(\Delta)=\min \{ m_i(\Delta_{[n]-x}) : x\in[n] \}$ for all $1\leq i < n-d$.
\item If $\Delta$ is 2-CM, then $m_{n-d}(\Delta)=M_{n-d}(\Delta)=n$.
\end{enumerate}
\end{lemma}
\proof Parts 1 and 2 are immediate from equations (\ref{M-interpr}) and (\ref{m-interpr}). Part 3 follows from \cite[Proposition~III.3.2(e)]{St96}.
\endproof

We are now in a position to prove Theorem \ref{matroid-multiplicity-thm}.

\smallskip\noindent {\it Proof of Theorem \ref{matroid-multiplicity-thm}: \ } The proof is by induction on $n$. The assertion clearly holds if $n=1$ or $n=2$. If $n>2$, then by Lemma \ref{matroid-properties} either $\Delta$ is a cone with apex $x$ or $\tilde{H}_{d-1}(\Delta; \field) \neq 0$.
In the former case, $\Delta_{[n]-x}$ is a $(d-2)$-dimensional matroid complex on $n-1$ vertices. Since in this case $f_{d-1}(\Delta)=f_{d-2}(\Delta_{[n]-x})$,
and $$
   M_i(\Delta)=M_i(\Delta_{[n]-x}) \mbox{ and } m_i(\Delta)=m_i(\Delta_{[n]-x}) \
\forall   1\leq i\leq n-d,
$$
all parts of the theorem follow from the induction hypothesis on $\Delta_{[n]-x}$.
In the latter case, each of $\Delta_{[n]-x}$ is a $(d-1)$-dimensional matroid complex on $n-1$ vertices. Thus we have \begin{eqnarray*}
f_{d-1}(\Delta) &=& \frac{1}{n-d}\sum_{x\in[n]} f_{d-1}(\Delta_{[n]-x}) \leq \frac{1}{n-d} \cdot \sum_{x\in[n]} \frac{\prod_{i=1}^{n-1-d} M_i(\Delta_{[n]-x})}{(n-1-d)!} \\ &\leq& \frac{1}{(n-d)!} \cdot n \cdot
      \prod_{i=1}^{n-1-d} \max_{x \in [n]} M_i(\Delta_{[n]-x}) = \frac{\prod_{i=1}^{n-d} M_i(\Delta)}{(n-d)!}.
\end{eqnarray*}
In the above equation, the first step is implied by the fact that every $(d-1)$-dimensional face has $d$ vertices, and hence is contained in exactly $n-d$ of the complexes $\Delta_{[n]-x}$. The second step is by induction hypothesis on $\Delta_{[n]-x}$, and the last step is an application of Parts 1 and 2 of Lemma \ref{M-and-m-properties}. Moreover, if equality
$f_{d-1}(\Delta) = (\prod_{i=1}^{n-d} M_i(\Delta))/(n-d)!$ holds, then all the inequalities in the above equation are equalities. Hence $$ f_{d-1}(\Delta_{[n]-x})=\frac{\prod_{i=1}^{n-d-1}
M_i(\Delta_{[n]-x})}{(n-d-1)!}
\mbox{ and } M_i(\Delta_{[n]-x})= M_i(\Delta) \quad \forall x\in [n], 1\leq i \leq n-d-1.
$$
The induction hypothesis on $\Delta_{[n]-x}$  then yields $M_i(\Delta)= M_i(\Delta_{[n]-x})=m_i(\Delta_{[n]-x})$ for all $x\in[n]$, and so, by Part 2 of Lemma \ref{M-and-m-properties}, $M_i(\Delta)=m_i(\Delta)$ for $i\leq n-d-1$. In addition, $M_{n-d}(\Delta)=m_{n-d}(\Delta)$ as $\Delta$ is 2-CM. Thus,  $\field[\Delta]$ has a pure resolution.

The proof of the lower bound (together with the treatment of equality) is
completely analogous and is omitted.
\endproof

We close this section with several remarks.
\begin{remark} {\rm Translating the circuit axiom for matroids into commutative
algebra leads to the following algebraic characterization:
\begin{itemize}
\item A proper squarefree monomial ideal $I$ is the face ideal of a matroid
complex if and only if for every pair of monomials $\mu_1, \mu_2\in I$ and for
every $i$ such that
$x_i$ divides both $\mu_1$ and $\mu_2$, the monomial $\lcm(\mu_1,\mu_2)/x_i$ is
in $I$ as well.
\end{itemize}
We refer to such an ideal as a {\bf matroid ideal}. In the commutative algebra 
literature the term ``matroid ideal" is sometimes used
for the ideal whose minimal generators correspond to the bases of a matroid.
Our notion of a matroid ideal is different as in our case the ideal is 
generated by
the circuits of the matroid rather than by its bases.

The above definition of a matroid ideal is
reminiscent of that of a squarefree strongly stable ideal (that is, the face
ideal of a shifted complex). In fact,  matroid complexes and  shifted complexes
share
certain properties. For instance, an induced subcomplex of a shifted complex is
shifted (cf, Lemma \ref{matroid-properties}(1)),
and  if $\Delta$ has a vanishing top homology and is a shifted  complex
on $[n]$ with respect to the ordering
$1\succ 2 \succ \cdots \succ n$, then $\Delta$ is a subcomplex of the cone
over $\Delta_{[n-1]}$  with apex $n$ (cf, Lemma \ref{matroid-properties}(3)).
Thus, the same reasoning as above
provides a new simple proof of the upper-bound part of the multiplicity
conjecture  for squarefree strongly stable ideals, the result
originally proved in~\cite{HerzSr98}.
}
\end{remark}

\begin{remark}      \label{0-homology}
{\rm
The same argument as in the proof of Theorem \ref{matroid-multiplicity-thm}
shows
that if the \mubc \ holds for all complexes of dimension $<d-1$ and for all
$(d-1)$-dimensional complexes with vanishing top homology, then it  holds
for all $(d-1)$-dimensional complexes.}
\end{remark}

\begin{remark} {\rm In view of Theorem \ref{matroid-multiplicity-thm} it is 
natural to ask
for which matroid complexes $\Delta$, does $\field[\Delta]$ have a pure
resolution.  Using either the methods of Theorem 3.11 of \cite{NPS}, or
Lemmas \ref{matroid-properties}, \ref{M-and-m-properties},
(\ref{M-interpr}), (\ref{m-interpr}), and the basic matroid properties of
duality, it is not hard to show that $k[\Delta]$ has a pure resolution if
and only if $\Delta$ is the matroid dual of a perfect matroid design.
Perfect matroid designs are matroids such that for every $r$ the
cardinality of a flat of rank $r$ is some fixed number $f(r).$  Thus, face
rings of duals of affine geometries and projective geometries provide a
large collection of examples of squarefree monomial ideals with  a pure
resolution.
  }
\end{remark}

\section{Lower bounds and the connectivity sequence}
As noted in the introduction, the key to the lower bounds are Theorems
\ref{main_thm} and \ref{q-estimate}. Before proving these, we recall an
alternative characterization of Cohen-Macaulay complexes due to Hochster.

\begin{theorem} \cite{Hochs} \label{alternative CM}
Let $\Delta$ be $(d-1)$-dimensional.  Then $\Delta$ is CM if and only if for
every subset $W$ of vertices
$$\tilde{H}_p(\Delta_W ; \field)=0, \mbox{ for all } p+(n-|W|) < d-1.$$
Thus, $\Skel_j(\Delta)$ is $q$-CM if and only if $$
\tilde{H}_p(\Delta_W ; \field)=0, \mbox{ for all } |W|-p \ge (n-q+1)-j+1.$$
\end{theorem}

We also need the following fact.
\begin{proposition} If $\Delta'$ is a $j$-dimensional CM complex, then
$\Skel_{j-1}(\Delta')$
is 2-CM. Hence if $\Delta'$ is $q$-CM complex, then its codimension one
skeleton is
$(q+1)$-CM.
\end{proposition}
\proof
By a result of Hibi \cite{Hibi88}, the codimension one skeleton of a CM complex
is level. In addition, $\Skel_{j-1}(\Delta')$ has a nonvanishing top homology,
as every $j$-face of $\Delta'$ is attached to a $(j-1)$-cycle
of $\Skel_{j-1}(\Delta')$. Thus
the assertion on 2-CM follows from the last paragraph of \cite[p.~94]{St96}.
For the second part apply the first one to $\Delta'$ with $q-1$ points removed.
\endproof

\begin{remark}       \label{Floystad}
A  corollary of the above is the previously mentioned result of Fl{\o}ystad
\cite{Floyst}, that for Cohen-Macaulay $\Delta$ with $\dim \Delta = d-1,$
$q_{d-1} < q_{d-2} < \dots < q_1 < q_0.$
\end{remark}

We are now ready to prove Theorems \ref{main_thm} and \ref{q-estimate}.
    For convenience, we repeat their statements.

\begin{numtheorem}{\ref{main_thm}} Let $\Delta$ be a  CM complex on $n$
vertices,
and let $(q_0, \ldots, q_{d-1})$ be the CM-connectivity sequence of $\Delta$,
where $d-1$ is the dimension of $\Delta$.
Then $$[n]-\{m_1(\Delta), \ldots, m_{n-d}(\Delta) \} = \{n-q_0+1, n-q_1+1,
\ldots, n-q_{d-1}+1\},
$$ and hence
\begin{equation} \label{lowerbound}
\frac{\prod_{i=1}^{n-d} m_i(\Delta)}{(n-d)!}
=\frac{n!}{(n-d)!} \cdot \frac{\prod_{i=1}^{n-d} m_i(\Delta)}{n!}=
\frac{n(n-1)\cdots(n-(d-1))}{(n-q_0+1)\cdots (n-q_{d-1}+1)}. \end{equation}
\end{numtheorem}

\proof
The $m$-sequence is a strictly increasing sequence of length $n-d$ of integers
contained in $[1,n].$  Hence, there are $d$ numbers skipped which we denote by
$s_0<s_1<\dots<s_{d-1}.$  We must prove that $q_j = n-s_j+1.$  We argue by
induction on $j$.
For $j=0$ this follows immediately from the fact that $m_1 \ge 2$ and $q_0=n.$

Let $m^\prime_i = m_i-i-1.$  So, $m^\prime_i$ is the dimension in which
$\tilde{H}_{|W_i| - i -1}(\Delta_{W_i};\field)$ is nonzero, where $W_i$ is a
subset of vertices of cardinality $m_i.$ Since the $m$-sequence is strictly
increasing, the $m^\prime$-sequence
   is nondecreasing. Define $t_j$ to be the largest $i$ such that $m^\prime_i <
j.$ With this definition, $s_j= t_j +j+ 1.$ (See Example \ref{example1} below
on how these invariants relate.)

Note that since $m'_i\geq j$ for all $i\geq t_j+1$, there can be no subsets $W$
of the
vertex set with $|W|\geq(t_j+1)+j=s_j$ and
$\tilde{H}_{j-1}(\Delta_W;\field)\neq 0$
(for all $j=0,1,\ldots,d-1$). Thus for a fixed $j$,
there is no subset $W$ with $|W|\geq s_j$ and
$\tilde{H}_{j-1}(\Delta_W;\field)\neq 0$,
no subset $W$ with $|W|\geq s_j-1 (\geq s_{j-1})$ and
$\tilde{H}_{j-2}(\Delta_W;\field)\neq 0$, etc. Theorem \ref{alternative CM}
then implies that $\Skel_j(\Delta)$ is $(n-s_j+1)$-CM, that is, $q_j\geq
n-s_j+1$.

Now if $j-1$ appears as some $m^\prime_i,$ then  $j-1=m'_{t_j}$,
$m_{t_j}=s_j-1$, and
$\tilde{H}_{j-1}(\Delta_{W_{t_j}})\neq 0$. Hence we also have $q_j\leq
n-|W_{t_j}|=n-s_j+1$, and so  $q_j=n-s_j+1$.
What happens if $j-1$ does not appear in the $m^\prime$-sequence?  In this case
$s_j=s_{j-1}+1$, and we infer from Remark \ref{Floystad}
   and induction hypothesis that
$q_j\leq q_{j-1}-1=(n-s_{j-1}+1)-1=n-s_j+1$ which again yields $q_j=n-s_j+1$.
\endproof

\begin{example}  \label{example1}
In this example $n=19$ and $d=9.$ $$\begin{array}{llllllllllll}
i & & 1 & 2&3&4&5&6&7&8&9&10\\
m_i& &2&3&4&6&7&11&13&16&17&18\\
m^\prime_i& &0&0&0&1&1&4&5&7&7&7
\end{array}$$
$$\begin{array}{lclclcl}
j& &t_j& &s_j& &q_j\\
   & & & & & &\\
0& &0& &1&  &19\\
1& &3& &5& &15\\
2& &5& &8& &12\\
3& &5& & 9& &11\\
4& &5& &10& &10\\
5& &6& &12& &8\\
6& &7& &14& &6\\
7& &7& &15& &5\\
8& &10& &19& &1
\end{array}$$
\end{example}

\begin{numtheorem}{\ref{q-estimate}} Let $\Delta$ be a  CM $(d-1)$-dimensional
complex,
and let $(q_0, \ldots, q_{d-1})$ be its CM-connectivity sequence.
Then
$$
e(\field[\Delta])=f_{d-1}(\Delta) \geq \frac{q_0q_1 \cdots q_{d-1}}{d!}.
$$

\end{numtheorem}

\proof
The proof is by induction on $d,$ with the initial case $d=0$ being
self-evident. Since links of $q$-CM-complexes are $q$-CM \cite{Baclaw} and the
$i$-skeleton of the link of a vertex $v$ is the link of $v$ in
the $(i+1)$-skeleton of $\Delta, q_i(\lk v) \ge q_{i+1}(\Delta).$ Therefore the
induction hypothesis implies that $$f_{d-2}(\lk v) \ge \frac{q_1 \cdots
q_{d-1}}{(d-1)!}.$$
Summing up over all $n = q_0$ vertices finishes the proof.
\endproof

Another theorem we will make a frequent use of is the following result.
Its first part is \cite[Theorem 1.2]{HerzSr98} and its second part is
\cite[Theorem 3.4]{HerzZheng}

\begin{theorem}     \label{quasi-pure}
   If $\field[\Delta]$ is Cohen-Macaulay and has a quasi-pure resolution, then
$\field[\Delta]$ satisfies the multiplicity conjecture.
Furthermore, if the multiplicity of $\field[\Delta]$ equals the lower bound or
the upper bound of the multiplicity
conjecture, then $\field[\Delta]$ has a pure resolution. \end{theorem}

With Theorems \ref{main_thm}, \ref{q-estimate}, and \ref{quasi-pure} in hand
we are ready to discuss 2-CM and $d$-CM complexes. We remark that although the
$M$-sequence  need not be strictly increasing in general,
it does strictly increase if $\Delta$ is CM \cite[Proposition~1.1]{BrunsHibi}.
Thus in the CM case we can also talk about skips in the $M$-sequence.

\begin{theorem}  \label{2CM-and-dCM}
Let $\Delta$ be a $(d-1)$-dimensional complex.
   If $\Delta$ is 2-CM with $q_1\leq n-d+1$ or $\Delta$ is $d$-CM, then
$\field[\Delta]$ satisfies the \mlbc.  Furthermore,
   if the bound is achieved,
then $\field[\Delta]$  has a pure resolution. \end{theorem}
\proof By virtue of Theorems \ref{main_thm} and \ref{q-estimate}, to prove the
\mlbc \ it suffices to verify that
$\frac{n(n-1)\cdots(n-(d-1))}{(n-q_0+1)\cdots (n-q_{d-1}+1)}
\leq \frac{q_0q_1 \cdots q_{d-1}}{d!},$
or equivalently (via $q_0=n$) that
\begin{equation}  \label{leq}
d! \prod_{i=1}^{d-1}(n-i) \leq \prod_{i=1}^{d-1} q_i(n-q_i+1).
\end{equation}
Assume first that $\Delta$ is 2-CM with $q_1\leq n-d+1$.
Then $n-d+1 \geq q_1 > q_2> \cdots > q_{d-1} \geq 2$ which implies that
$$
q_{d-i} \in [i+1, n-d+1] \subseteq [i+1, n-i] \quad \forall 1\leq i \leq d-1.
$$
Since $f(q):=q(n-q+1)$ is a concave function of $q$ symmetric about
$q=(n+1)/2$, it follows that $q_{d-i}(n-q_{d-i}+1)\geq (i+1)(n-i)$ for all
$1\leq i \leq d-1$.
Multiplying these inequalities up over all $1\leq i \leq d-1$ yields
(\ref{leq}).

When can the bound $f_{d-1}=(\prod_{i=1}^{n-d} m_i)/(n-d)!$ be achieved?
Since for all $i<d-1$, $q_{d-i} \in [i+1, n-d+1] \subset [i+1, n-i)$,
analysis of the above inequalities and the proof of Theorem \ref{q-estimate}
reveals that this can happen only if $q_1\in \{d, n-d+1\}$, $q_{d-i}=i+1$ for
$i< d-1$,
and $f_{d-3}(\lk E)=(q_2\cdots q_{d-1})/(d-2)!=d-1$ for every edge $E$.
As $\lk E$ is a $(d-3)$-dimensional 2-CM complex, the latter condition implies
that
for every edge $E$, $\lk E$ is the boundary of a $(d-2)$-simplex, which in turn
implies that $\Delta$
itself is the boundary of a $d$-simplex as long as $d>3$, and hence that
$\field[\Delta]$ has a pure resolution. This completes the treatment of
equality in the $d>3$ case. If $d=2$ then the resolution is always quasi-pure
and the result follows from
   Theorem \ref{quasi-pure}.
What if $d=3$? Then either (i) $q_2=2$, $q_1=3$, and $f_1(\lk v)=2\cdot 3/2!=3$
for every vertex $v$
or (ii) $q_2=2$ and $q_1=n-2$. In the former case, the link of each vertex must
be the boundary of a 2-simplex, and hence $\Delta$ itself must be
the boundary of a 3-simplex.
In the latter case the $1^{st}$ skip in the $m$-sequence is $s_1=3$.
Since the $2^{nd}$ (and the last) skip in the $M$-sequence is at least
3, the resolution is quasi-pure, and we are done by Theorem \ref{quasi-pure}.

We now turn to the case of $d$-CM $\Delta$. In this case $d\leq q_{d-1}< \cdots
< q_2 <q_1 \leq n-1$. Thus
$q_i\in[d,n-i] \subseteq [i+1, n-i]$, and the same computation as in the 2-CM
case implies (\ref{leq}). Moreover, if the multiplicity
lower bound is achieved, then $q_{d-1}\in\{d, n-d+1\}$ and $q_i=n-i$ for all
$i<d-1$.
The latter implies (via Theorem \ref{main_thm}) that all integers from 1 to
$d-1$ are skipped from the $m$-sequence, and hence that the resolution is
quasi-pure. \endproof

The hypothesis $q_1 \le n-d+1$  is very mild as it only requires that
there be some $(d-1)$-subset of vertices which is disconnected.   For
instance, (reduced) order complexes of all of the following posets satisfy
this condition: face posets of 2-CM cell complexes with the intersection 
property
(that is, intersection of any two faces is a face; this class includes face 
posets of all polytopes and face posets of all 2-CM simplicial complexes), 
geometric lattices,
supersolvable lattices with nonzero M\"obius function on every interval, rank 
selected subposets of any of these.
We remark that the \mubc \  for the order complexes of
face posets of all simplicial
complexes was very recently verified by Kubitzke and Welker~\cite{KubWelk}.

\section{One- and two-dimensional complexes}
The goal of this section is to establish the multiplicity conjecture
for one- and two-dimensional complexes. We start with the \mubc.
This requires the following strengthening of Theorem \ref{quasi-pure}.
\begin{definition}
   The face ring $\field[\Delta]$ of a $(d-1)$-dimensional complex $\Delta$ is
{\bf almost Cohen-Macaulay} if the length of its  minimal free resolution is at
most $n-d+1.$ (Equivalently, $\field[\Delta]$ is almost CM if the codimension
one skeleton of $\Delta$ is CM.)
\end{definition}

\begin{theorem} \cite[Theorem 1.5]{HerzSr98}            \label{quasi-pure'}
    If $\field[\Delta]$ is almost Cohen-Macaulay and has a quasi-pure 
resolution,
then $\field[\Delta]$ satisfies the \mubc.
   \end{theorem}

\begin{theorem}
If $\Delta$ is one or two-dimensional, then  $\field[\Delta]$ satisfies the
\mubc.  Furthermore, if $\Delta$ is one-dimensional and

$$f_1 = \frac{1}{(n-2)!} \prod^{n-2}_{i=1} M_i,$$
then $\field[\Delta]$ is Cohen-Macaulay and has a pure resolution.
\end{theorem}

\proof
First we consider $\dim \Delta = 1.$ Using (\ref{M-interpr}) and
(\ref{m-interpr}) we see that in this case $m_i, M_i \in \{i+1, i+2\}$, and so
$\field[\Delta]$ has a quasi-pure resolution.
Also since $M_i\leq n$, it follows that $i\leq n-1$. Thus
$\field[\Delta]$ is almost Cohen-Macaulay, and hence satisfies the \mubc.
   If $\Delta$ is not connected and has $t$ components, let $\Delta^\prime$ be
any connected complex obtained by adding $t-1$ edges to $\Delta.$  Then
$\Delta^\prime$ is connected, has the same 1-cycles as $\Delta$ and has more
edges than $\Delta.$ Since the $M_i$ only depend on the cardinality of the
cycles, this reasoning shows that if equality occurs, then $\Delta$ is
connected and hence Cohen-Macaulay.  In this case, Theorem \ref{quasi-pure}
implies that $\field[\Delta]$ has a pure resolution.

Now assume that $\Delta$ is $2$-dimensional, connected, and
$\tilde{H}_2(\Delta; \field) = 0$.  Again, by (\ref{M-interpr}) and
(\ref{m-interpr}),  $\Delta$ is almost Cohen-Macaulay and has a quasi-pure
resolution, and hence satisfies the \mubc. When $\Delta$ is not connected we
can add edges as above to get a connected complex with the same number of
triangles and identical $M_i.$ What if $\tilde{H}_2 \neq 0?$ Remark
\ref{0-homology} completes the proof in this case.
\endproof

For one-dimensional complexes, i.e. graphs,  the $M_i$ encode the size of the
smallest circuit.  If the graph is acyclic, then the multiplicity conjecture
gives the best possible bound, $f_1 \le n-1.$ When the graph contains a
triangle Conjecture \ref{multiplicity-conj} says that $f_1 \le n(n-1)/2$ which,
in view of the complete graph, is best possible.
   However, for all other possible smallest circuit sizes the asymptotic upper
bound for $f_1$ is known to be much less than that given by the
   multiplicity upper bound conjecture.  Determination of the optimal upper
bounds with this data is an area of ongoing research.
   Triangle free graphs have at most $n^2/4$ edges (this is Mantel's theorem
\cite[p.~30]{vanLint} --- a special case of Tur{\'a}n's theorem). The maximum
number of edges in a graph without triangles or squares is asymptotically
bounded above by $n \sqrt{n-1}/2$ \cite{GarKwongLaz}.

We now turn to the \mlbc.
   Recall that for a CM complex $\Delta$, $h_i(\Delta)\geq 0$ for all $i$
\cite{St75}, \cite[Cor.~II.3.2]{St96}.

\begin{theorem} \label{1-2-lower}
   If $\Delta$ is a $1$ or $2$-dimensional CM complex, then it satisfies the
\mlbc.  Furthermore,
   if
   \begin{equation} \label{2dim_pure}
   f_{d-1} = \frac{1}{(n-d)!} \prod^{n-d}_{i=1} m_i,
   \end{equation}
then $\field[\Delta]$  has a pure resolution. \end{theorem}

    \proof
    If $\Delta$ is one-dimensional, then there are only two skips in the $m$ and
$M$-sequence and $1$ is the $0^{th}$ skip in both. Therefore, the minimal
resolution is quasi-pure and Theorem~\ref{quasi-pure}
applies. So, we assume that $\dim \Delta=2.$
   There are several cases to consider.

    \begin{enumerate}
    \item \label{case1}
      $\tilde{H}_2(\Delta; \field) = 0.$  Under these conditions there are only
two skips in the $M$-sequence (as $M_{n-3}<n$), and so $\field[\Delta]$ has a
quasi-pure resolution.

      \item
      $q_1 = n-1.$  This is equivalent to the 1-skeleton being the complete 
graph
on $n$ vertices.  Now  the $1^{st}$ skip of the $m$-sequence is $2.$ Once again
$\field[\Delta]$ has a quasi-pure resolution.

      \item
      $q_2 \ge 2.$  If $q_1 = n-1,$ then the previous case holds. Otherwise,
$q_1\leq n-2$ and  Theorem \ref{2CM-and-dCM} applies.

      \item $q_2 = 1$ and $q_1$ is $2$ or $3.$ When $q_1=2$ Theorem
\ref{main_thm} says that we must show that $f_2 \ge n-2.$  However, for any
two-dimensional CM
complex, $h_0 = 1$, $h_1 = n-3$, and $h_2, h_3\geq 0$. Hence $f_2 = h_0 + h_1 +
h_2 + h_3 \ge n-2.$   If $q_1 = 3,$ then we must show that $f_2 \ge n-1.$ If
$f_2 < n-1,$ then the $h$-vector must be $(1,n-3,0,0).$ This implies that
$\dim_\field \tilde{H}_2(\Delta;\field) = h_3 = 0$ and that case \ref{case1}
applies.

      \item $q_2=1$ and $q_1 \ge 4.$ By Theorem \ref{q-estimate} applied to the
$1$-skeleton, $f_1 \ge n q_1/2.$ Using (\ref{h-vector}), $h_2 \ge n q_1/2 -
2(n-3) -3$ and hence,
      $$f_2 = 1 + (n-3) + h_2 +h_3 \ge \frac{n (q_1-2)+2}{2}.$$
   Comparing this to the \mlbc \  via  (\ref{lowerbound}),
   $$\frac{n (q_1-2)+2}{2} - \frac{(n-1)(n-2)}{n-q_1+1} =
\frac{[n-q_1-1][(q_1-4)n +2]}{2(n-q_1+1)} \ge 0.$$
    \end{enumerate}
   The last inequality holds since $q_1 \ge 4$ and any Cohen-Macaulay complex
$\Delta$ with $\dim \Delta \ge 1$ has at least $q_1 +1$ vertices. If
(\ref{2dim_pure}) holds, then $n=q_1+1$ and the $1$-skeleton is the complete
graph.  In addition, we must have $0=h_3 (=\dim\tilde{H}_3(\Delta; \field))$
since this is implicit in the above estimate. Thus the skips for both the $m$
and the $M$-sequences are $1,2$ and $n.$ \endproof

\section{Gorenstein complexes}
In this section we discuss the multiplicity conjecture for
Gorenstein complexes. Since every such complex is
the join of a simplex and a Gorenstein* complex, it suffices to treat the case
of Gorenstein* complexes only.
As in the previous section we start with the \mubc.
This will require the following facts.
\begin{lemma}          \label{prod-M-expression}
Let $\Delta$ be a $(d-1)$-dimensional Gorenstein* complex on $n$ vertices,
and let $(q_0, \ldots, q_{d-1})$ be its CM-connectivity sequence. If $d\geq 3$,
then
\begin{enumerate}
\item $q_{d-1}=2$ and $q_{d-2}\leq 5$.
\item The sequence $M_1, \ldots M_{n-d}$ is strictly increasing and satisfies
$[n]-\{M_1, \ldots, M_{n-d} \} = \{q_{d-1}-1, q_{d-2}-1, \ldots, q_1-1,
q_0-1\}.
$ Hence
$$
\frac{\prod_{i=1}^{n-d} M_i}{(n-d)!}
=\frac{n!}{(n-d)!} \cdot \frac{\prod_{i=1}^{n-d} M_i}{n!}=
n \cdot \frac{(n-2)\cdots(n-(d-1))}{(q_{1}-1) \cdots (q_{d-2}-1)}. $$
\end{enumerate}
\end{lemma}
\proof If $\Delta$ is a 2-dimensional Gorenstein* complex on $n$ vertices, then
   $f_1(\Delta)= 3n-6$.
(This is a consequence of the Euler relation
$f_0-f_1+f_2=2$ and the fact that every 2-face has exactly three edges, while
every edge is contained in exactly
two 2-faces.) Thus $\Delta$ has a vertex of degree $\leq 5$, and hence the
1-skeleton of $\Delta$ is at most 5-connected. To see that
   $q_{d-2}\leq 5$ for a general $d>3$, note that the link of a
$(d-4)$-dimensional face $F$ of $\Delta$ is a 2-dimensional Gorenstein*
complex and that $q_{d-2}(\Delta) \leq q_1 (\lk F)$ \cite{Baclaw}.

Since $\Delta$ is Gorenstein*, the minimal free resolution of $I_\Delta$ is
self-dual.
Thus $$M_i+m_{n-d-i}=n   \quad \mbox{for all $0\leq i \leq n-d$},
$$
   where we set $m_0=M_0=0$.
This fact together with Theorem \ref{main_thm} implies Part 2.
Finally, since $2\leq M_1 < \ldots < M_{n-d}$, we have
$1\in [n]-\{M_1, \ldots M_{n-d}\}$,
and so $q_{d-1}-1=1$.
\endproof

In view of the lemma, we refer to $q_{d-i-1}-1$ ($0\leq i \leq d-1$) as the
{\bf $i$-th skip} in the $M$-sequence. Note that for a Gorenstein* complex,
$q_{0}-1=n-q_{d-1}+1$, and so the  $(d-1)$-th skip of the $M$-sequence
coincides with that of the $m$-sequence.
Note also that $q_{d-2}>3$ implies that $2$ is not skipped from the
$M$-sequence, and hence that $M_1=2$.

\begin{theorem}   \label{3-4-Gor-upper}
If $\Delta$ is a three- or four-dimensional Gorenstein* complex, then
$\field[\Delta]$ satisfies the \mubc. Furthermore, if the bound is achieved
then $\field[\Delta]$  has a pure resolution.
\end{theorem}

\proof We can assume without loss of generality that $\field[\Delta]$ does not
have a quasi-pure resolution (since otherwise the result follows from Theorem
\ref{quasi-pure}).

First consider the case $d-1=\dim \Delta =4$. Since the resolution is not
quasi-pure, either
the first skip in the $m$-sequence occurs after the second skip
in the $M$-sequence or the second skip in $m$ occurs after the third skip in
$M$, that is, either $n-q_1+1>q_2-1$ or $n-q_2+1>q_1-1$.
Both of these expressions are equivalent to $(q_1-1)+(q_2-1)\leq n-1$ which
together with $q_1\geq q_2+1$  yields
$$
(q_1-1)(q_2-1)\leq \frac{n}{2} \left(\frac{n}{2}-1\right).
$$
Substituting this inequality along with $q_3\leq 5$ in Part 2
of Lemma \ref{prod-M-expression}, we obtain
$$
\frac{\prod_{i=1}^{n-5} M_i}{(n-5)!} =
\frac{n(n-2)(n-3)(n-4)}{(q_1-1)(q_2-1)(q_3-1)} \geq (n-3)(n-4) =2 {n-3 \choose
2} \geq f_4,$$
where the last inequality is  the Upper Bound Theorem for Gorenstein*
complexes \cite{St75}.
Furthermore, equality $(\prod_{i=1}^{n-5} M_i)/(n-5)!=f_4$ is impossible under
above assumptions as it would imply
$q_3-1=4$, and  hence $M_1=2$, on one hand, and $f_4=2{n-3 \choose 2}$, and
hence 2-neighborliness of the Gorenstein* complex, on the other.

Now assume that $\Delta$ is 3-dimensional. Then
the assumption that $\field[\Delta]$
does not have a quasi-pure resolution implies that $q_1-1<n-q_1+1$, and so
$2(q_1-1)\leq n-1$.
We treat separately two cases: $q_2=3$ and $q_2\in\{4,5\}$.
If $q_2=3$, then we have $$
\frac{\prod_{i=1}^{n-4} M_i}{(n-4)!} =
\frac{n(n-2)(n-3)}{(q_1-1)(q_2-1)} \geq \frac{n(n-2)(n-3)}{n-1} >
(n-2)(n-3) > \frac{n(n-3)}{2}
\geq f_3,$$
where the first step is Lemma \ref{prod-M-expression}
and the last step is again the Upper Bound Theorem.
The $q_2\in\{4,5\}$ case requires a little bit more work.
First we note that in this case $n-1 \geq 2(q_1-1)\geq 2q_2 \geq 8$, and so
$n\geq 9$.
Also since $q_2> 3$, we have $M_1=2$.
Thus $\Delta$ is a 3-dimensional
{\em flag} complex. In particular, the 1-skeleton of $\Delta$ is a
5-clique-free graph.
Tur{\'a}n's theorem \cite[Theorem~4.1]{vanLint} then implies that
$$
f_1(\Delta) \leq \frac{n(n-\lfloor n/4 \rfloor)}{2} \leq \frac{n(n-2)}{2}.
$$
It remains to note that for a 3-dimensional Gorenstein* complex
$f_3=f_1-n$. (This is a consequence of the Euler relation
$f_0-f_1+f_2-f_3=0$ and the fact that every 3-face has exactly four 2-faces,
while every 2-face is contained in exactly
two 3-faces). Hence
$$
f_3 =f_1-n \leq \frac{n(n-4)}{2} <  \frac{n(n-2)(n-3)}{2(n-1)}
\leq \frac{n(n-2)(n-3)}{(q_1-1)(q_2-1)}=\frac{\prod_{i=1}^{n-4} M_i}{(n-4)!}.
$$
\endproof

Generalizing case $d-1=3$, $q_2=3$ of the above proof leads to the following
observation.

\begin{proposition}
The \mubc \ holds for the following classes of Gorenstein* complexes.
\begin{enumerate}
\item $(d-1)$-dimensional complexes with $M_1\geq \lfloor d/2\rfloor+1$;
\item $(d-1)$-dimensional complexes on $n$ vertices, where $n-d=4$ and $d\leq
22$.
\end{enumerate}
\end{proposition}
\proof To prove the first assertion note that $M_{i}\geq m_{i}=n-M_{n-d-i}$,
and so $M_i+M_{n-d-i}\geq n$ for $i\geq 1$. A routine
computation using these inequalities
together with $\lfloor d/2\rfloor+1 \leq M_1 < M_2 < \cdots < M_{n-d}=n$ shows
that
the upper bound of the multiplicity conjecture is at least as large as the
bound provided by the Upper Bound Theorem for Gorenstein* complexes \cite{St75}
implying the result.

Similarly, for the second assertion one verifies that either $\field[\Delta]$
has a quasi-pure resolution or the upper bound of the multiplicity conjecture
is at least as large as the bound of the Upper Bound Theorem. We omit the
details.
\endproof

Finally we discuss the lower bounds in the three- and four-dimensional cases.
\begin{theorem}
   If $\Delta$ is a $3$ or $4$-dimensional Gorenstein* complex, then it 
satisfies
the \mlbc.  Furthermore,
   if the lower bound is achieved
then $\field[\Delta]$  has a pure resolution. \end{theorem}

\proof $\Delta$ is Gorenstein*, hence 2-CM, and so as long as
$q_1\leq n- \dim\Delta$ Theorem \ref{2CM-and-dCM} applies.
Thus we assume for the rest of the proof that $q_1\geq n-\dim\Delta+1$.

By Theorem \ref{quasi-pure}, we can also assume that $\field[\Delta]$ does not
have a quasi-pure resolution which, as was observed
in the proof of Theorem \ref{3-4-Gor-upper}, is equivalent to
\begin{eqnarray}
    \nonumber
   2(q_1-1)\leq n-1 \quad &\mbox{ if }& \dim\Delta=3, \mbox{ and }\\
   (q_1-1)+(q_2-1) \leq n-1 \quad &\mbox{ if }& \dim\Delta=4.   \label{4}
\end{eqnarray}
This finishes the proof of the 3-dimensional case as the inequalities $q_1\geq
n-2$ and $2(q_1-1)\leq n-1$ imply
$2(n-3)\leq 2(q_1-1)\leq n-1$, and hence $n\leq 5$, on one hand, and
$n-1\geq 2(q_1-1)\geq2(q_3+1)=6$ on the other.

Assume now $\dim\Delta=4$. If $q_1\geq n-2$, then, by (\ref{4}), $q_2\leq 3$
which contradicts the fact that $q_2\geq q_4+2=4$. Thus $q_1=n-3$, and $q_2=4$,
$q_3=3$, $q_2=2$, and via Theorem \ref{main_thm} to complete
the proof we must show that $f_4>n(n-4)/4$.

By Theorem \ref{q-estimate} applied to the 1-skeleton,
$f_1\geq nq_1/2=n(n-3)/2$. Hence, by (\ref{h-vector}),
$h_2=f_1-4f_0+10 \geq (n^2-11n+20)/2$. The Dehn-Sommerville relations
\cite{Klee} then yield
$$
f_4=2(h_0+h_1+h_2)\geq 2[1+(n-5)+ (n^2-11n+20)/2]=n^2-9n+12.$$
This completes the proof if $n\geq 9$, since in this case
$n^2-9n+12>n(n-4)/4$.

What if $n\leq 8$? Since $n-3=q_1\geq q_2+1=5$, $n$ must equal 8. Then
$h_0=1$, $h_1=8-5=3$, and $h_2>0$. Thus
$
f_4=2(h_0+h_1+h_2)\geq 2[1+3+1]=10 > 8(8-4)/4,
$
as desired.
\endproof

\end{document}